\input amstex
\input amsppt.sty
\magnification=\magstep1
\hsize=30truecc
\vsize=22.2truecm
\baselineskip=16truept
\NoBlackBoxes
\TagsOnRight \pageno=1 \nologo
\def\Z{\Bbb Z}
\def\N{\Bbb N}

\def\Q{\Bbb Q}

\def\l{\left}
\def\r{\right}
\def\bg{\bigg}
\def\({\bg(}
\def\[{\bg\lfloor}
\def\){\bg)}
\def\]{\bg\rfloor}
\def\t{\text}
\def\f{\frac}

\def\p{\ (\roman{mod}\ p)}

\def\bi{\binom}
\def\eq{\equiv}

\def\ls{\leqslant}
\def\gs{\geqslant}
\def\mo{\roman{mod}}

\def\Proof{\noindent{\it Proof}}

\def\Remark{\medskip\noindent{\it  Remark}}

\def\Ack{\medskip\noindent {\bf Acknowledgment}}
\hbox {J. Number Theory 131(2011), no.\,11, 2219-2238.}\bigskip
\topmatter
\title On congruences related to central binomial coefficients\endtitle
\author Zhi-Wei Sun\endauthor
\leftheadtext{Zhi-Wei Sun} \rightheadtext{Congruences related to
central binomial coefficients}
\affil Department of Mathematics, Nanjing University\\
 Nanjing 210093, People's Republic of China
  \\  zwsun\@nju.edu.cn
  \\ {\tt http://math.nju.edu.cn/$\sim$zwsun}
\endaffil
\abstract
It is known that
$$\sum_{k=0}^\infty\f{\bi{2k}k}{(2k+1)4^k}=\f{\pi}2\ \ \t{and}\ \ \sum_{k=0}^\infty\f{\bi{2k}k}{(2k+1)16^k}=\f{\pi}3.$$
In this paper we obtain their $p$-adic analogues such as
$$\sum_{p/2<k<p}\f{\bi{2k}k}{(2k+1)4^k}\eq3\sum_{p/2<k<p}\f{\bi{2k}k}{(2k+1)16^k}\eq pE_{p-3}\ (\mo\ p^2),$$
where $p>3$ is a prime and $E_0,E_1,E_2,\ldots$ are Euler numbers.
Besides these, we also deduce some other congruences related to central binomial coefficients.
In addition, we pose some conjectures one of which states that for any odd prime $p$ we have
$$\sum_{k=0}^{p-1}\bi{2k}k^3\eq\cases4x^2-2p\ (\mo\ p^2)&\t{if}\ (\f p7)=1\ \&\ p=x^2+7y^2\ (x,y\in\Z),
\\0\ (\mo\ p^2)&\t{if}\ (\f p7)=-1,\ \t{i.e.,}\ p\eq3,5,6\ (\mo\ 7).\endcases$$
\endabstract
\thanks 2010 {\it Mathematics Subject Classification}. Primary 11B65;
Secondary 05A10, 11A07, 11B68, 11E25.
\newline\indent Supported by the National Natural Science
Foundation (grant 10871087) and the Overseas Cooperation Fund (grant 10928101) of China.
\endthanks
\endtopmatter
\document

\heading{1. Introduction}\endheading

The following three series related to $\pi$ are well known (cf. [Ma]):
$$\sum_{k=0}^\infty\f{\bi{2k}k}{(2k+1)4^k}=\f{\pi}2,\ \ \ \ \sum_{k=0}^\infty\f{\bi{2k}k}{(2k+1)16^k}=\f{\pi}3,$$
and
$$\sum_{k=0}^\infty\f{\bi{2k}k}{(2k+1)^2(-16)^k}=\f{\pi^2}{10}.$$
These three identities can be easily shown by using $1/(2k+1)=\int_0^1 x^{2k}dx$.
In March 2010 the author [Su2] suggested that
$$\sum_{k=0}^\infty\f{\bi{2k}k}{(2k+1)^3 16^k}=\f{7\pi^3}{216}$$
via a public message to Number Theory List,
and then Olivier Gerard pointed out there is a computer proof via certain math. softwares
like {\tt Mathematica} (version 7).
Our main goal in this paper is to investigate  $p$-adic analogues of the above identities for powers of $\pi$.

For a prime $p$ and an integer $a\not\eq0\ (\mo\ p)$,
we let $q_p(a)$ denote the Fermat quotient $(a^{p-1}-1)/p$.
For an odd prime $p$ and an integer $a$, by $(\f ap)$ we mean the Legendre symbol.
As usual, harmonic numbers refer to those $H_n=\sum_{0<k\ls n}1/k$ with $n\in\N=\{0,1,2,\ldots\}$.
Recall that Euler numbers $E_0,E_1,E_2,\ldots$ are integers defined by
$E_0=1$ and the recursion:
$$ \sum^n\Sb k=0\\2\mid k\endSb \bi nk E_{n-k}=0\ \ \ \t{for}\ n=1,2,3,\ldots.$$
And Bernoulli numbers $B_0,B_1,B_2,\ldots$ are rational numbers given by $B_0=1$ and
$$\sum^n_{k=0}\bi {n+1}k B_k=0\ \ \ \ (n=1,2,3,\ldots).$$

Now we state our first theorem which gives certain $p$-adic analogues of the first and the second identities
mentioned at the beginning of this section.

\proclaim{Theorem 1.1} Let $p$ be an odd prime.

{\rm (i)} We have
$$\sum_{k=0}^{(p-3)/2}\f{\bi{2k}k}{(2k+1)4^k}\eq(-1)^{(p+1)/2}q_p(2)\ (\mo\ p^2),\tag1.1$$
and
$$\sum_{p/2<k<p}\f{\bi{2k}k}{(2k+1)4^k}\eq pE_{p-3}\ (\mo\ p^2)\tag1.2$$
which is equivalent to the congruence
$$\sum_{k=1}^{(p-1)/2}\f{4^k}{(2k-1)\bi{2k}k}\eq E_{p-3}+(-1)^{(p-1)/2}-1\ (\mo\ p).\tag1.3$$

{\rm (ii)} Suppose $p>3$. Then
$$\sum_{k=0}^{(p-3)/2}\f{\bi{2k}k}{(2k+1)16^k}\eq0\ (\mo\ p^2),\tag1.4$$
and
$$\sum_{p/2<k<p}\f{\bi{2k}k}{(2k+1)16^k}\eq \f p3E_{p-3}\ (\mo\ p^2)\tag1.5$$
which is equivalent to the congruence
$$\sum_{k=1}^{(p-1)/2}\f{16^k}{k(2k-1)\bi{2k}k}\eq\f 83 E_{p-3}\ (\mo\ p).\tag1.6$$
\endproclaim

\Remark\ 1.1. Motivated by the work of H. Pan and Z. W. Sun [PS], and Sun and R. Tauraso [ST1, ST2], the author [Su1]
 managed to determine
$\sum_{k=0}^{p^a-1}\bi{2k}k/m^k$
modulo $p^2$, where $p$ is a prime, $a$ is a positive integer, and $m$ is any integer not divisible by $p$.
See also [SSZ], [G-Z] and [Su3] for related results on $p$-adic valuations.

\medskip

The congruences in Theorem 1.1 are somewhat sophisticated.
Now we deduce some easier congruences via combinatorial identities.
Using the software {\tt Sigma}, we find the identities
$$\sum_{k=0}^n\bi nk\f{(-1)^k}{(2k+1)^2}=\f{4^n}{(2n+1)\bi{2n}n}\sum_{k=0}^n\f1{2k+1},$$
$$\sum_{k=0}^n\f{(-1)^k}{(k+1)\bi nk}=n+1-(n+1)\sum_{k=1}^n\f{1-2(-1)^k}{(k+1)^2},$$
and
$$\align n\sum_{k=2}^n\f{(-1)^k}{(k-1)^2\bi nk}=&\sum_{k=2}^n\f{1-2k+(-1)^k(1-k+2k^2)}{k(k-1)^2}
\\=&\f{1+(-1)^n}n-\sum_{k=1}^{n-1}\f{1+2(-1)^k}{k^2}.
\endalign$$
If $p=2n+1$ is an odd prime, then
$$\bi nk\eq\bi{-1/2}k=\f{\bi{2k}k}{(-4)^k}\ (\mo\ p)\quad\t{for all}\ k=0,\ldots,p-1.$$
Thus, from the above three identities we deduce for any prime $p>3$ the congruences
$$\sum_{k=0}^{(p-3)/2}\f{\bi{2k}k}{(2k+1)^24^k}\eq(-1)^{(p+1)/2}\ \f{q_p(2)^2}2\ (\mo\ p),\tag1.7$$
$$\sum_{k=2}^{(p-1)/2}\f{4^k}{(k-1)^2\bi{2k}k}\eq 8E_{p-3}-4-12\l(\f{-1}p\r)\ (\mo\ p)\tag1.8$$
and
$$\sum_{k=0}^{(p-1)/2}\f{4^k}{(k+1)\bi{2k}k}\eq\l(\f{-1}p\r)(4-2E_{p-3})-2\ (\mo\ p).\tag1.9$$
Note that the series $\sum_{k=0}^\infty 4^k/((k+1)\bi{2k}k)$ diverges while {\tt Mathematica} (version 7) yields
$$\sum_{k=0}^\infty\f{\bi{2k}k}{(2k+1)^24^k}=\f{\pi}4\log2
\ \ \ \t{and}\ \ \ \sum_{k=2}^\infty \f{4^k}{(k-1)^2\bi{2k}k}=\pi^2-4$$
the latter of which appeared in [Sp].

Let $p$ be an odd prime. By a known result (see, e.g., [I]),
$$\sum_{k=0}^{p-1}\f{\bi{2k}k^3}{64^k}\eq a(p)\ (\mo\ p^2),$$
where the sequence $\{a(n)\}_{n\gs1}$ is defined by
$$\sum_{n=1}^\infty a(n)q^n=q\prod_{n=1}^\infty(1-q^{4n})^6.$$
Clearly, $a(p)=0$ if $p\eq 3\ (\mo\ 4)$.

Recall that Catalan numbers are those integers
$$C_k=\f1{k+1}\bi{2k}k=\bi{2k}k-\bi{2k}{k+1}\ \ (k=0,1,2,\ldots).$$
They have many combinatorial interpretations (see, e.g., [St2, pp.\,219-229]).

Now we present our second theorem.
\proclaim{Theorem 1.2} Let $p$ be an odd prime.

{\rm (i)} We have
$$\sum_{k=0}^{p-1}\f{k^3\bi{2k}k^3}{64^k}\eq\cases0\ (\mo\ p)&\t{if}\ p\eq1\ (\mo\ 4),
\\-\f1{640}(\f{p+1}4!)^{-4}\ (\mo\ p)&\t{if}\ p\eq3\ (\mo\ 4).\endcases\tag1.10$$
If $p>3$ and $p\eq 3\ (\mo\ 4)$, then
$$\sum_{k=0}^{p-1}\bi{p-1}k\f{\bi{2k}k^3}{(-64)^k}\eq0\ (\mo\ p^2).\tag1.11$$

{\rm (ii)} We have
$$\sum_{k=0}^{p-1}\f{C_k^2}{16^k}\eq-3\ (\mo\ p),\tag1.12$$
and
$$\sum_{k=0}^{p-1}\f{C_k^3}{64^k}\eq\cases7\ (\mo\ p)&\t{if}\ p\eq1\ (\mo\ 4),
\\7-\f32(\f{p+1}4!)^{-4}\ (\mo\ p)&\t{if}\ p\eq3\ (\mo\ 4).\endcases\tag1.13$$
Also,
$$\sum_{k=0}^{p-1}\f{\bi{2k}kC_k}{32^k}\eq\cases p\ (\mo\ p^2)&\t{if}\ p\eq1\ (\mo\ 4),
\\p+(4p+2^{p}-6)\binom{(p-3)/2}{(p-3)/4}\ (\mo\ p^2)&\t{if}\ p\eq3\ (\mo\ 4).\endcases\tag1.14$$
\endproclaim
\Remark\ 1.2. Let $p$ be an odd prime. We conjecture that if $p\eq1\ (\mo\ 4)$ and $p>5$ then
$$\sum_{k=0}^{p^a-1}\f{k^3\bi{2k}k^3}{64^k}\eq0\ (\mo\ p^{2a})\quad \ \t{for all}\ a=1,2,3,\ldots.$$
We also conjecture that $\sum_{k=0}^{(p-1)/2}kC_k^3/16^k\eq2p-2\pmod{p^2}$ if $p\eq1\pmod3$, and
$\sum_{k=0}^{(p-1)/2}C_k^3/64^k\eq8\pmod{p^2}$ if $p\eq1\pmod4$.
\medskip

In the next section we are going to provide several lemmas.
Theorems 1.1 and 1.2 will be proved in Sections 3 and 4 respectively.
Section 5 contains some open conjectures of the author for further research.

\heading{2. Some lemmas}\endheading

For $n\in\N$ the Chebyshev polynomial $U_n(x)$ of the second kind is given by
$$U_n(\cos\theta)=\f{\sin((n+1)\theta)}{\sin\theta}.$$
It is well known that
$$U_n(x)=\sum_{k=0}^{\lfloor n/2\rfloor}\bi{n-k}k(-1)^k(2x)^{n-2k}.$$

\proclaim{Lemma 2.1} For $n\in\N$, we have the identities
$$\sum_{k=0}^n\bi{n+k}{2k}\f{(-4)^k}{2k+1}=\f{(-1)^n}{2n+1}\tag 2.1$$
and
$$\sum_{k=0}^n\bi{n+k}{2k}\f{(-1)^k}{2k+1}=\cases(-1)^n/(2n+1)&\t{if}\ 3\nmid 2n+1,
\\2(-1)^{n-1}/(2n+1)&\t{if}\ 3\mid 2n+1.\endcases\tag2.2$$
\endproclaim
\Proof.  Note that
$$U_{2n}(x)=\sum_{k=0}^n\bi{2n-k}{2n-2k}(-1)^k(2x)^{2n-2k}=\sum_{j=0}^n\bi{n+j}{2j}(-1)^{n-j}(2x)^{2j}.$$
Thus
$$\align\sum_{k=0}^n\bi{n+k}{2k}\f{(-4)^k}{2k+1}=&\int_0^1\sum_{k=0}^n\bi{n+k}{2k}(-4)^kx^{2k}dx
=(-1)^n\int_0^1U_{2n}(x)dx
\\=&(-1)^n\int_{\pi/2}^0U_{2n}(\cos\theta)(-\sin\theta)d\theta
\\=&(-1)^n\int^{\pi/2}_0\sin((2n+1)\theta)d\theta
\\=&\f{-(-1)^n}{2n+1}\cos((2n+1)\theta)\bigg|_0^{\pi/2}=\f{(-1)^n}{2n+1}.
\endalign$$
Similarly,
$$\align\sum_{k=0}^n\bi{n+k}{2k}\f{(-1)^k}{2k+1}=&\int_0^1\sum_{k=0}^n\bi{n+k}{2k}(-1)^kx^{2k}dx
=(-1)^n\int_0^1U_{2n}\l(\f x2\r)dx
\\=&(-1)^n\int_{\pi/2}^{\pi/3}U_{2n}(\cos\theta)(-2\sin\theta)d\theta
\\=&-2(-1)^n\int_{\pi/2}^{\pi/3}\sin((2n+1)\theta)d\theta
\\=&\f{2(-1)^n}{2n+1}\cos((2n+1)\theta)\bigg|_{\pi/2}^{\pi/3}=\f{2(-1)^n}{2n+1}\cos\l(\f{2n+1}3\pi\r)
\\=&\cases(-1)^n/(2n+1)&\t{if}\ 3\nmid 2n+1,
\\2(-1)^{n-1}/(2n+1)&\t{if}\ 3\mid 2n+1.\endcases
\endalign$$
This concludes the proof. \qed

\proclaim{Lemma 2.2} Let $p=2n+1$ be an odd prime. For $k=0,\ldots,n$ we have
$$\bi{n+k}{2k}\eq\f{\bi{2k}k}{(-16)^k}\ (\mo\ p^2).\tag2.3$$
\endproclaim
\Proof. As observed by the author's brother Z. H. Sun,
$$\align\bi{n+k}{2k}=&\f{\prod_{0<j\ls k}(p^2-(2j-1)^2)}{4^k(2k)!}
\\\eq&\f{\prod_{0<j\ls k}(-(2j-1)^2)}{4^k(2k)!}
=\f{\bi{2k}k}{(-16)^k}\ (\mo\ p^2).
\endalign$$
We are done. \qed
\Remark\ 2.1. Using Lemma 2.2 and the identity
$$\sum_{k=0}^n\f{\bi{n+k}{2k}(-2)^k}{2k+1}=\f{(1+i)(-i)^n(1+(-1)^{n-1}i)}{2(2n+1)},$$
we can deduce for any prime $p>3$ that
$$\sum_{k=0}^{(p-3)/2}\f{\bi{2k}k}{(2k+1)8^k}
\eq-\l(\f{-2}p\r)\f{q_p(2)}2+\l(\f{-2}p\r)\f p8q_p^2(2)\ (\mo\ p^2).$$

\proclaim{Lemma 2.3} Let $p$ be any odd prime. Then
$$\sum_{k=1}^{(p-1)/2}\f{4^k}{k\bi{2k}k}\eq 2\l((-1)^{(p-1)/2}-1\r)\ (\mo\ p).\tag2.4$$
\endproclaim
\Proof. Clearly (2.4) holds for $p=3$.

Now assume that $p>3$. We can even show a stronger congruence
$$\f12\sum_{k=1}^{(p-1)/2}\f{4^k}{k\bi{2k}k}\eq (-1)^{(p-1)/2}(1-p\,q_p(2)+p^2q_p(2)^2)-1\ (\mo\ p^3).$$
Let us employ a known identity (cf. [G, (2.9)])
$$\sum_{k=1}^n\f{2^{2k-1}}{k\bi{2k}k}=\f{2^{2n}}{\bi{2n}n}-1$$
which can be easily proved by induction. Taking $n=(p-1)/2$ and noting that
$$(-1)^n\bi{2n}n\eq 4^{p-1}\ (\mo\ p^3)$$
by Morley's congruence ([Mo]), we get
$$\align\f12\sum_{k=1}^{(p-1)/2}\f{4^k}{k\bi{2k}k}\eq &\f{(-1)^{(p-1)/2}}{1+p\,q_p(2)}-1
\\\eq&(-1)^{(p-1)/2}(1-p\,q_p(2)+p^2q_p(2)^2)-1\ (\mo\ p^3).
\endalign$$
This ends the proof. \qed

\proclaim{Lemma 2.4} For any $n\in\N$, we have the identity
$$\sum_{k=-n}^n\f{(-1)^k}{(2k+1)^2}\bi{2n}{n+k}=\f{16^n}{(2n+1)^2\bi{2n}n}.\tag2.5$$
\endproclaim
\Proof. Let $u_n$ and $v_n$ denote the left-hand side and the right-hand side of (2.5)
respectively. By the well-known Zeilberger algorithm (cf. [PWZ]),
$$(2n+3)(2n+5)^2u_{n+2}-16(n+2)(2n+3)^2u_{n+1}
+64(n+1)(n+2)(2n+1)u_n=0$$
for all $n=0,1,2,\ldots$. It is easy to verify that $\{v_n\}_{n\gs0}$ also satisfies this recurrence.
Since $u_0=v_0=1$ and $u_1=v_1=8/9$, by the recursion we have $u_n=v_n$ for all $n\in\N$.
 \qed

\Remark\ 2.2. (2.5) was discovered by the author during his study of Delannoy numbers (cf. [Su5]).
The reader may consult [GZ] and [ZG] for some other combinatorial identities obtained via solving recurrence relations.

\proclaim{Lemma 2.5} For any $n\in\N$ we have
$$\sum_{k=0}^n\bi{2n-k}k(-1)^k=\l(\f{1-n}3\r)\tag2.6$$
and
$$\sum_{k=0}^n\bi{2n-k}k\f1{(-4)^{k}}=\f{2n+1}{4^n}.\tag2.7$$
\endproclaim
\Remark\ 2.3. (2.6) and (2.7) are known identities, see (1.75) and (1.73) of [G].

\proclaim{Lemma 2.6} Let $p>3$ be a prime. Then
$$\sum_{0<k\ls\lfloor p/6\rfloor}\f{(-1)^k}{k^2}\eq (-1)^{(p-1)/2}10E_{p-3}\ (\mo\ p).\tag2.8$$
\endproclaim
\Proof. Recall that the Euler polynomial of degree $n$ is defined by
$$E_n(x)=\sum_{k=0}^n\bi nk\f{E_k}{2^k}\l(x-\f12\r)^{n-k}.$$
It is well known that
$$E_n(1-x)=(-1)^nE_n(x),\ \ E_n(x)+E_n(x+1)=2x^n,$$
and $$ E_n(x)=\f2{n+1}\l(B_{n+1}(x)-2^{n+1}B_{n+1}\l(\f x2\r)\r),$$
where $B_m(x)$ denotes the Bernoulli polynomial of degree $m$.

Note that $E_{p-3}(0)=\f2{p-2}(1-2^{p-2})B_{p-2}=0$ and $E_{p-3}(5/6)=E_{p-3}(1/6)$. Thus
$$\align2\sum_{0<k\ls\lfloor p/6\rfloor}\f{(-1)^k}{k^2}\eq&\sum_{k=0}^{\lfloor p/6\rfloor}(-1)^k(2k^{p-3})
\\=&\sum_{k=0}^{\lfloor p/6\rfloor}\l((-1)^kE_{p-3}(k)-(-1)^{k+1}E_{p-3}(k+1)\r)
\\=&E_{p-3}(0)-(-1)^{\lfloor p/6\rfloor+1}E_{p-3}\l(\l\lfloor \f p6\r\rfloor+1\r)
\\\eq&(-1)^{\lfloor p/6\rfloor}E_{p-3}\l(\f16\r)\ (\mo\ p).
\endalign$$
Evidently $\lfloor p/6\rfloor\eq(p-1)/2\ (\mo\ 2)$. As
$E_n(1/6)=2^{-n-1}(1+3^{-n})E_n$ for all $n=0,2,4,\ldots$ (see, e.g., G. J. Fox [F]), we have
$$E_{p-3}\l(\f 16\r)=2^{2-p}(1+3^{3-p})E_{p-3}\eq 2(1+3^2)E_{p-3}=20E_{p-3}\ (\mo\ p).$$
Therefore (2.8) follows from the above. \qed

\heading{3. Proof of Theorem 1.1}\endheading

\medskip
\noindent{\it Proof of Theorem 1.1}.  (a) Set $n=(p-1)/2$. By Lemmas 2.1 and 2.2,
$$\align\sum_{k=0}^{n-1}\f{\bi{2k}k}{(2k+1)4^k}\eq&\sum_{k=0}^{n-1}\bi{n+k}{2k}\f{(-4)^k}{2k+1}=\f{(-1)^n-(-4)^n}{2n+1}
\\=&(-1)^{n}\f{1-2^{p-1}}p=(-1)^{n+1}q_p(2)\ (\mo\ p^2).
\endalign$$
This proves (1.1). When $p=2n+1>3$, again by Lemmas 2.1 and 2.2, we have
$$\sum_{k=0}^{n-1}\f{\bi{2k}k}{(2k+1)16^k}\eq\sum_{k=0}^{n-1}\bi{n+k}{2k}\f{(-1)^k}{2k+1}=0\ (\mo\ p^2)$$
and hence (1.4) holds.

(b) For $k\in\{1,\ldots,(p-1)/2\}$, it is clear that
$$\align\f1p\bi{2(p-k)}{p-k}=&\f1p\times\f{p!\prod_{s=1}^{p-2k}(p+s)}{((p-1)!/\prod_{0<t<k}(p-t))^2}
\\\eq&\f{(k-1)!^2}{(p-1)!/(p-2k)!}\eq-\f{(k-1)!^2}{(2k-1)!}=-\f2{k\bi{2k}k}\ (\mo\ p).
\endalign$$
Therefore
$$\align &\f1p\sum_{p/2<k<p}\f{\bi{2k}k}{(2k+1)4^k}=\sum_{k=1}^{(p-1)/2}\f{\bi{2(p-k)}{p-k}/p}{(2(p-k)+1)4^{p-k}}
\\\eq&-2\sum_{k=1}^{(p-1)/2}\f{4^{k-1}}{(1-2k)k\bi{2k}k}=\f12\sum_{k=1}^{(p-1)/2}\f{4^k}{k(2k-1)\bi{2k}k}\ (\mo\ p).
\endalign$$
Similarly,
$$\f1p\sum_{p/2<k<p}\f{\bi{2k}k}{(2k+1)16^k}\eq\f18\sum_{k=1}^{(p-1)/2}\f{16^k}{k(2k-1)\bi{2k}k}\ (\mo\ p)$$
and hence (1.5) and (1.6) are equivalent.
Observe that
$$\sum_{k=1}^{(p-1)/2}\f{4^k}{k(2k-1)\bi{2k}k}
=2\sum_{k=1}^{(p-1)/2}\f{4^k}{(2k-1)\bi{2k}k}-\sum_{k=1}^{(p-1)/2}\f{4^k}{k\bi{2k}k}.$$
Thus, in view of (2.4), both (1.2) and (1.3) are equivalent to the congruence
$$\sum_{k=1}^{(p-1)/2}\f{4^k}{k(2k-1)\bi{2k}k}\eq 2E_{p-3}\ (\mo\ p)\tag3.1$$
which holds trivially when $p=3$.

(c) Now we prove (3.1) for $p>3$. It is easy to see that
$$(n+1)(2(n+1)-1)\bi{2(n+1)}{n+1}=2(2n+1)^2\bi{2n}n$$
for any $n\in \N$. Thus
$$\sum_{k=1}^{(p-1)/2}\f{4^k}{k(2k-1)\bi{2k}k}=\sum_{k=0}^{(p-3)/2}\f{4^{k+1}}{2(2k+1)^2\bi{2k}k}.$$
In view of Lemma 2.4,
$$\align\sum_{n=0}^{(p-3)/2}\f{4^n}{(2n+1)^2\bi{2n}n}=&\sum_{n=0}^{(p-3)/2}\f1{4^n}\sum_{k=-n}^n\f{(-1)^k}{(2k+1)^2}\bi{2n}{n-k}
\\=&\sum_{k=-(p-3)/2}^{(p-3)/2}\f{(-1)^k}{(2k+1)^24^{|k|}}\sum_{n=|k|}^{(p-3)/2}\f{\bi{2n}{n-|k|}}{4^{n-|k|}}.
\endalign$$
For $k\in\{0,\ldots,(p-3)/2\}$, with the help of Lemma 2.5 we have
$$\align \sum_{n=k}^{(p-3)/2}\f{\bi{2n}{n-k}}{4^{n-k}}=&\sum_{r=0}^{(p-3)/2-k}\f{\bi{2k+2r}r}{4^r}
=\sum_{r=0}^{(p-3)/2-k}\f{\bi{-2k-r-1}r}{(-4)^r}
\\\eq&\sum_{r=0}^{(p-1)/2-k}\f{\bi{p-1-2k-r}r}{(-4)^r}-\f1{(-4)^{(p-1)/2-k}}
\\=&\f{p-2k-(-1)^{(p-1)/2-k}}{4^{(p-1)/2-k}}\eq\f{(-1)^{(p+1)/2-k}-2k}{4^{-k}}\ (\mo\ p).
\endalign$$
Therefore
$$\align\sum_{n=0}^{(p-3)/2}\f{4^n}{(2n+1)^2\bi{2n}n}\eq&\sum_{k=-(p-3)/2}^{(p-3)/2}\f{(-1)^k}{(2k+1)^2}\l((-1)^{(p+1)/2-k}-2|k|\r)
\\=&\sum_{k=0}^{(p-3)/2}\f{(-1)^k}{(2k+1)^2}\l((-1)^{(p+1)/2-k}-2k\r)
\\&+\sum_{k=1}^{(p-3)/2}\f{(-1)^{-k}}{(-2k+1)^2}\l((-1)^{(p+1)/2+k}-2k\r)
\\\eq&\sum_{k=1}^{(p-1)/2}\f{(-1)^{k-1}}{(2k-1)^2}\l((-1)^{(p+1)/2-k+1}-2(k-1)\r)
\\&+\sum_{k=1}^{(p-1)/2}\f{(-1)^{k}}{(2k-1)^2}\l((-1)^{(p+1)/2+k}-2k\r)\ \ (\mo\ p)
\endalign$$
and hence
$$\align\sum_{n=0}^{(p-3)/2}\f{4^n}{(2n+1)^2\bi{2n}n}\eq&\sum_{k=1}^{(p-1)/2}\f{(-1)^k}{(2k-1)^2}\l(2(-1)^{(p+1)/2-k}+2(k-1)-2k\r)
\\=&4(-1)^{(p+1)/2}\sum\Sb 1\ls k\ls(p-1)/2\\k\eq(p-1)/2\ (\mo\ 2)\endSb
\f1{(2k-1)^2}\ \ (\mo\ p).
\endalign$$
Since $p>3$ and $\sum_{k=1}^{p-1}1/(2k)^2\eq\sum_{k=1}^{p-1}1/k^2\ (\mo\ p)$, we have
$$2\sum_{k=1}^{(p-1)/2}\f1{k^2}\eq\sum_{k=1}^{(p-1)/2}\l(\f1{k^2}+\f1{(p-k)^2}\r)=\sum_{k=1}^{p-1}\f1{k^2}\eq0\ (\mo\ p)$$
and hence
$$\align \sum\Sb 1\ls k\ls(p-1)/2\\2k+1\eq p\ (\mo\ 4)\endSb\f1{(2k-1)^2}
\eq&\sum\Sb 1\ls k\ls(p-1)/2\\p+1-2k\eq2\ (\mo\ 4)\endSb\f1{(p+1-2k)^2}
=\sum^{p-1}\Sb k=1\\k\eq2\ (\mo\ 4)\endSb\f1{k^2}
\\\eq&-\sum^{p-1}\Sb k=1\\4\mid k\endSb\f1{k^2}=-\f1{16}\sum_{k=1}^{\lfloor p/4\rfloor}\f1{k^2}
\ \ (\mo\ p).
\endalign$$
As $\sum_{k=1}^{\lfloor p/4\rfloor}1/k^2\eq(-1)^{(p-1)/2}4E_{p-3}\ (\mo\ p)$ by Lehmer [L, (20)], from the above we obtain that
$$\sum_{n=0}^{(p-3)/2}\f{4^n}{(2n+1)^2\bi{2n}n}\eq 4(-1)^{(p+1)/2}\f{(-1)^{(p-1)/2}4E_{p-3}}{-16}=E_{p-3}\ \ (\mo\ p)$$
and hence (3.1) holds.

(d) Finally we show (1.6) for $p>3$. In view of Lemmas 2.4 and 2.5, arguing as in (c) we get
$$\align &\f18\sum_{k=1}^{(p-1)/2}\f{16^k}{k(2k-1)\bi{2k}k}=\sum_{n=0}^{(p-3)/2}\f{16^n}{(2n+1)^2\bi{2n}n}
\\\eq&\sum_{k=-(p-3)/2}^{(p-3)/2}\f{(-1)^k}{(2k+1)^2}\sum_{r=0}^{(p-3)/2-|k|}\bi{p-1-2|k|-r}r(-1)^r
\\=&\sum_{k=-(p-3)/2}^{(p-3)/2}\f{(-1)^k}{(2k+1)^2}\(\l(\f{|k|-(p-3)/2}3\r)-(-1)^{(p-1)/2-|k|}\)
\\=&\sum_{k=-(p-3)/2}^{(p-3)/2}\f{(-1)^k}{(2k+1)^2}\(\l(\f{p-2|k|}3\r)+(-1)^{(p+1)/2-|k|}\)\pmod{p}.
\endalign$$
Observe that
$$\align \sum_{k=-(p-3)/2}^{(p-3)/2}\f1{(2k+1)^2}=&\sum_{k=0}^{(p-3)/2}\f1{(2k+1)^2}+\sum_{k=1}^{(p-3)/2}\f1{(-2k+1)^2}
\\=&2\sum_{k=1}^{(p-1)/2}\f1{(2k-1)^2}-\f1{(p-2)^2}
\\\eq&2\(\sum_{k=1}^{p-1}\f1{k^2}-\sum_{k=1}^{(p-1)/2}\f1{(2k)^2}\)-\f14\eq-\f14\ \ (\mo\ p)
\endalign$$
and
$$\align &\sum_{k=-(p-3)/2}^{(p-3)/2}\f{(-1)^k}{(2k+1)^2}\l(\f{p-2|k|}3\r)
\\=&\sum_{k=0}^{(p-3)/2}\f{(-1)^k}{(2k+1)^2}\l(\f{p+k}3\r)+\sum_{k=1}^{(p-3)/2}\f{(-1)^k}{(-2k+1)^2}\l(\f{p+k}3\r)
\\=&\sum_{k=1}^{(p-1)/2}\f{(-1)^k}{(2k-1)^2}\l(\l(\f{p+k}3\r)-\l(\f{p+k-1}3\r)\r)-\f{(-1)^{(p-1)/2}}{(p-2)^2}\l(\f{p+(p-1)/2}3\r)
\\\eq&\sum_{k=1}^{(p-1)/2}\f{(-1)^k}{(2k-1)^2}-3\sum^{(p-1)/2}
\Sb k=1\\3\mid p+k+1\endSb\f{(-1)^k}{(2k-1)^2}+\f{(-1)^{(p+1)/2}}4\ \ (\mo\ p).
\endalign$$
Therefore
$$\align&\f18\sum_{k=1}^{(p-1)/2}\f{16^k}{k(2k-1)\bi{2k}k}
\\\eq&\sum_{k=1}^{(p-1)/2}\f{(-1)^k}{(p-(2k-1))^2}
-3\sum^{(p-1)/2}\Sb k=1\\3\mid 2k-1-p\endSb\f{(-1)^k}{(p-(2k-1))^2}
\\=&\sum_{k=1}^{(p-1)/2}\f{(-1)^{(p+1)/2-k}}{(2k)^2}-3\sum_{0<k\ls\lfloor p/6\rfloor}\f{(-1)^{(p+1)/2-3k}}{(6k)^2}
\pmod{p}.
\endalign$$
Since
$$\sum_{k=1}^{(p-1)/2}\f{(-1)^k}{k^2}\eq\sum_{k=1}^{(p-1)/2}\f{(-1)^k+1}{k^2}=\f12\sum_{k=1}^{\lfloor p/4\rfloor}\f1{k^2}
\eq 2(-1)^{(p-1)/2}E_{p-3}\ (\mo\ p),$$
with the help of Lemma 2.6 we finally get
$$\f18\sum_{k=1}^{(p-1)/2}\f{16^k}{k(2k-1)\bi{2k}k}\eq-\f{E_{p-3}}2+\f{10}{12}E_{p-3}=\f{E_{p-3}}3\ (\mo\ p)$$
which proves (1.6).
\medskip

Combining (a)-(d) we have completed the proof of Theorem 1.1. \qed

\heading{4. Proof of Theorem 1.2}\endheading

\proclaim{Lemma 4.1} For any $n\in\N$ we have
$$\sum_{k=0}^{2n}(-1)^k\bi{2n}k^3=(-1)^n\f{(3n)!}{(n!)^3}.\tag4.1$$
\endproclaim
\Proof. By Dixon's identity (cf. [St1, p.\,45]) we have
$$\sum_{k=-n}^n(-1)^{k}\bi{2n}{n+k}^3=\f{(3n)!}{(n!)^3},$$
which is equivalent to the desired identity. \qed

\proclaim{Lemma 4.2 {\rm ([DPSW, (2)])}} For any positive odd integer $n$ we have the identity
$$\sum_{k=0}^n\bi nk^3(-1)^kH_k=\f{(-1)^{(n+1)/2}}3\cdot\f{(3n)!!}{(n!!)^3},\tag4.2$$
where $(2m+1)!!$ refers to $\prod_{k=0}^m(2k+1)$.
\endproclaim

\proclaim{Lemma 4.3} For each $n=1,2,3,\ldots$, we have
$$\sum_{k=0}^n\bi{n+k}{2k}\f{C_k}{(-2)^k}=\cases(-1)^{(n-1)/2}C_{(n-1)/2}/2^n&\t{if}\ 2\nmid n,
\\0&\t{if}\ 2\mid n.\endcases\tag4.3$$
\endproclaim
\Proof. The desired identity can be easily proved by the WZ method (cf. [PWZ]); in fact,
if we denote by $S(n)$ the sum of the left-hand side or the right-hand side of (4.3), then
we have the recursion $S(n+2)=-nS(n)/(n+3)\ (n=1,2,3,\ldots)$. \qed

\medskip
\noindent{\it Proof of Theorem 1.2}. Let us recall that
$$\bi{(p-1)/2}k\eq\bi{-1/2}k=\f{\bi{2k}k}{(-4)^k}\ \ \t{for}\ k=0,1,\ldots,p-1.$$
Note also that for any positive odd integer $n$ we have
$$\sum_{k=0}^n(-1)^k\bi nk^3=\f12\sum_{k=0}^n\l((-1)^k\bi nk^3+(-1)^{n-k}\bi n{n-k}^3\r)=0.$$
These two basic facts will be frequently used in the proof.

(i) Clearly,
$$\align \sum_{k=0}^{p-1}\f{k^3\bi{2k}k^3}{64^k}\eq&\sum_{k=0}^{(p-1)/2}(-1)^kk^3\bi{(p-1)/2}k^3
\\=&-\l(\f{p-1}2\r)^3\sum_{k=1}^{(p-1)/2}(-1)^{k-1}\bi{(p-3)/2}{k-1}^3
\\\eq&\f18\sum_{k=0}^{(p-3)/2}(-1)^k\bi{(p-3)/2}k^3\ (\mo\ p).
\endalign$$
So, if $p\eq 1\ (\mo\ 4)$ then $(p-3)/2$ is odd and hence $\sum_{k=0}^{p-1}k^3\bi{2k}k^3/64^k\eq0\ (\mo\ p)$.
When $p=4n+3$ with $n\in\N$, applying Lemma 4.1 we get
$$\align 8\sum_{k=0}^{p-1}\f{k^3\bi{2k}k^3}{64^k}\eq&(-1)^n\f{(3n)!}{(n!)^3}
=\f{(-1)^{n}((p+1)/4)^3}{((p+1)/4)!^3}\times\f{(p-1)!}{\prod_{0<k<p-3n}(p-k)}
\\\eq&\f{(-1)^{n+1}}{64((p+1)/4)!^3(-1)^{p-1-3n}(p-1-3n)!}
\\\eq&-\f1{64((p+1)/4)!^4(p+5)/4}\ (\mo\ p).
\endalign$$
So (1.10) holds.

For $k=0,1,\ldots,p-1$, clearly
$$\bi{p-1}k(-1)^k=\prod_{0<j\ls k}\l(1-\f pj\r)\eq1-pH_k\ (\mo\ p^2).$$
When $p>3$ and $p\eq3\ (\mo\ 4)$, $\sum_{k=0}^{p-1}\bi{2k}k^3/64^k\eq0\ (\mo\ p^2)$
as mentioned in the first section, hence with the help of Lemma 4.2 we get
$$\align\sum_{k=0}^{p-1}\bi{p-1}k\f{\bi{2k}k^3}{(-64)^k}\eq&\sum_{k=0}^{p-1}(1-pH_k)\f{\bi{2k}k^3}{64^k}
\\\eq&-p\sum_{k=0}^{(p-1)/2}\bi{(p-1)/2}k^3(-1)^kH_k
\\\eq&-p\f{(-1)^{(p+1)/4}}3\times\f{(3(p-1)/2)!!}{((p-1)/2)!!^3}\eq0\ (\mo\ p^2).
\endalign$$
This proves (1.11) for $p\eq3\ (\mo\ 4)$ with $p\not=3$.

(ii) Below we set $n=(p-1)/2$ and want to show (1.12)-(1.14). Note that  $C_k\eq0\ (\mo\ p)$ when $n<k<p-1$. Also,
$$C_{p-1}=\f1p\bi{2p-2}{p-1}=\f1{2p-1}\bi{2p-1}p\eq-\prod_{k=1}^{p-1}\f{p+k}k\eq-1\ (\mo\ p).$$
Thus,
$$\sum_{k=0}^{p-1}\f{C_k^2}{16^k}\eq\sum_{k=0}^{n}\f{C_k^2}{16^k}+1\eq \sum_{k=0}^n\f1{(k+1)^2}\bi nk^2+1\ (\mo\ p)$$
and
$$\sum_{k=0}^{p-1}\f{C_k^3}{64^k}\eq\sum_{k=0}^{n}\f{C_k^3}{64^k}-1
\eq\sum_{k=0}^n\f{(-1)^k}{(k+1)^3}\bi nk^3-1\ (\mo\ p).$$
Clearly,
$$\align&(n+1)^2\sum_{k=0}^{n}\f1{(k+1)^2}\bi nk^2
\\=&\sum_{k=0}^n\bi{n+1}{k+1}^2=\sum_{k=0}^{n+1}\bi{n+1}k^2-1
\\=&\sum_{k=0}^{n+1}\bi{n+1}k\bi{n+1}{n+1-k}-1
\\=&\bi{2n+2}{n+1}-1\ (\t{by the Chu-Vandermonde identity (cf. [GKP, p.\,169])})
\\=&\bi{p+1}{(p+1)/2}-1=\f {2p}{(p-1)/2}\bi{p-1}{(p-3)/2}-1\eq-1\ (\mo\ p)
\endalign$$
and
$$\align&-(n+1)^3\sum_{k=0}^n\f{(-1)^k}{(k+1)^3}\bi nk^3
\\=&\sum_{k=0}^n(-1)^{k+1}\bi{n+1}{k+1}^3=\sum_{k=0}^{n+1}(-1)^k\bi{n+1}k^3-1.
\endalign$$

If $p\eq1\ (\mo\ 4)$, then $n+1$ is odd and hence
$$\sum_{k=0}^{n+1}(-1)^k\bi{n+1}k^3=0.$$
When $p=4m-1$ with $m\in\Z$, by Lemma 4.1
$$\sum_{k=0}^{n+1}(-1)^k\bi{n+1}k^3=\sum_{k=0}^{2m}(-1)^k\bi{2m}k^3=(-1)^m\f{(3m)!}{(m!)^3},$$
and in the case $m>1$ we have
$$\align(-1)^m(3m)!=&(-1)^m\f{(p-1)!}{\prod_{0<k<m-1}(p-k)}
\\\eq&-\f1{(m-2)!}=-\f{m(m-1)}{m!}\eq\f3{16(m!)}\ (\mo\ p).
\endalign$$
Therefore, if $p\eq3\ (\mo\ 4)$ then
$$\sum_{k=0}^{n+1}(-1)^k\bi{n+1}k^3\eq\f3{16}\l(\f{p+1}4!\r)^{-4}\ (\mo\ p).$$

By the above,
$$\sum_{k=0}^{p-1}\f{C_k^2}{16^k}\eq1-\f1{(n+1)^2}=1-\f 4{(p+1)^2}\eq-3\ (\mo\ p).$$
If $p\eq1\ (\mo\ 4)$, then
$$\sum_{k=0}^{p-1}\f{C_k^3}{64^k}\eq\f1{(n+1)^3}-1=\f 8{(p+1)^3}-1\eq7\ (\mo\ p).$$
If $p\eq3\ (\mo\ 4)$, then
$$\sum_{k=0}^{p-1}\f{C_k^3}{64^k}\eq\f{\f3{16}(\f{p+1}4!)^{-4}-1}{-(n+1)^3}-1\eq7-\f32\l(\f{p+1}4!\r)^{-4}\ (\mo\ p).$$
This proves (1.12) and (1.13).

With the help of Lemma 2.2, we have
$$\align\sum_{k=0}^{p-1}\f{\bi{2k}kC_k}{32^k}=&\f{pC_{p-1}^2}{32^{p-1}}+\sum_{k=0}^{p-2}\f{\bi{2k}kC_k}{32^k}
\\\eq& p+\sum_{k=0}^n\bi{n+k}{2k}\f{C_k}{(-2)^k}\ (\mo\ p^2).
\endalign$$

If $p\eq1\ (\mo\ 4)$, then $n=(p-1)/2$ is even and hence
$$\sum_{k=0}^n\bi{n+k}{2k}\f{C_k}{(-2)^k}=0$$
by Lemma 4.3.

Now assume that $p\eq3\ (\mo\ 4)$. In view of Lemma 4.3,
$$\align\sum_{k=0}^n\bi{n+k}{2k}\f{C_k}{(-2)^k}
=&(-1)^{(n-1)/2}\f{C_{(n-1)/2}}{2^n}
\\=&\f{(-1)^{(p-3)/4}}{2^{(p-1)/2}((p-3)/4+1)}\bi{(p-3)/2}{(p-3)/4}
\\\eq&\f{4(p-1)}{1+(\f2p)(2^{(p-1)/2}-(\f2p))}\bi{(p-3)/2}{(p-3)/4}\ (\mo\ p^2).
\endalign$$
Note that
$$\align&\f{4(p-1)}{1+(\f2p)(2^{(p-1)/2}-(\f2p))}
\\\eq&4(p-1)\l(1-\l(\f2p\r)\l(2^{(p-1)/2}-\l(\f2p\r)\r)\r)
\\\eq&(4p-4)\l(1-\f{2^{p-1}-1}2\r)\eq 4p-4+2(2^{p-1}-1)\ (\mo\ p^2).
\endalign$$

By the above, the congruence (1.14) also holds.
We are done. \qed

\heading{5. Some open conjectures}\endheading

In this section we pose some conjectures for further research.

Motivated by the identities $\sum_{k=0}^\infty\bi{2k}k/((2k+1)16^k)=\pi/3$,
$$\sum_{k=0}^\infty\f{\bi{2k}k}{(2k+1)^2(-16)^k}=\f{\pi^2}{10}\ \ \t{and}\ \
\sum_{k=0}^\infty\f{\bi{2k}k}{(2k+1)^316^k}=\f{7\pi^3}{216},
$$
we formulate the following conjecture based on our computation via {\tt Mathematica}.

 \proclaim{Conjecture 5.1} Let $p>5$ be a prime.
Then
$$\sum_{k=0}^{(p-3)/2}\f{\bi{2k}k}{(2k+1)16^k}\eq(-1)^{(p-1)/2}\l(\f{H_{p-1}}{12}+\f {3p^4}{160}B_{p-5}\r)\ (\mo\ p^5)$$
and
$$\sum_{k=0}^{(p-3)/2}\f{\bi{2k}k}{(2k+1)^3 16^k}\eq(-1)^{(p-1)/2}\l(\f{H_{p-1}}{4p^2}+\f{p^2}{36}B_{p-5}\r)\ (\mo\ p^3).$$
We also have
$$\align\sum_{k=0}^{(p-3)/2}\f{\bi{2k}k}{(2k+1)^2(-16)^k}\eq&\f{H_{p-1}}{5p}\ (\mo\ p^3),
\\\sum_{p/2<k<p}\f{\bi{2k}k}{(2k+1)^2(-16)^k}\eq&-\f p4B_{p-3}\ (\mo\ p^2).
\endalign$$
\endproclaim
\Remark\ 5.1. It is known that $H_{p-1}\eq-p^2B_{p-3}/3\ (\mo\ p^3)$ for any prime $p>3$ (see, e.g., [S]).
Thus the first congruence in the conjecture is a refinement of (1.4).

\medskip

Motivated by the known identities
$$\sum_{k=1}^\infty\f{2^k}{k^2\bi{2k}k}=\f{\pi^2}8\ \ \ \t{and}\ \ \ \sum_{k=1}^\infty\f{3^k}{k^2\bi{2k}k}=\f29\pi^2$$
(cf. [Ma]), we raise the following related conjecture.

\proclaim{Conjecture 5.2} Let $p$ be an odd prime. Then
$$\sum_{k=1}^{p-1}\f{\bi{2k}k}{k2^k}\eq-\f{H_{(p-1)/2}}2+\f7{16}p^2B_{p-3}\ (\mo\ p^3)$$
When $p>3$, we have
$$\align\sum_{k=1}^{p-1}\f{(-2)^k}{k^2}\bi{2k}k\eq&-2q_p(2)^2\pmod p,
\\p\sum_{k=1}^{p-1}\f{2^k}{k^2\bi{2k}k}\eq&-q_p(2)+\f{p^2}{16}B_{p-3}\ (\mo\ p^3),
\\\sum_{k=1}^{p-1}\f{\bi{2k}k}{k3^k}\eq&-2\sum^{p-1}\Sb k=1\\k\not\eq p\,(\mo\ 3)\endSb\f1k\ (\mo\ p^3),
\endalign$$
and
$$p\sum_{k=1}^{p-1}\f{3^k}{k^2\bi{2k}k}\eq-\f 32q_p(3)+\f 49p^2B_{p-3}\ (\mo\ p^3).$$
\endproclaim

Now we propose three more conjectures.

\proclaim{Conjecture 5.3}
Let $p$ be an odd prime. Then
$$\align&\sum_{k=0}^{p-1}\bi{2k}k^3
\\\eq&\cases4x^2-2p\ (\mo\ p^2)&\t{if}\ (\f p7)=1\ \&\ p=x^2+7y^2\ \t{with}\ x,y\in\Z,
\\0\ (\mo\ p^2)&\t{if}\ (\f p7)=-1,\ \t{i.e.,}\ p\eq3,5,6\ (\mo\ 7).\endcases
\endalign$$
\endproclaim

\Remark\ 5.2. Let $p$ be an odd prime with $(\f p7)=1$. As $(\f{-7}p)=1$, and the quadratic field $\Q(\sqrt{-7})$
has class number one, $p$ can be written uniquely in the form
$$\f{a+b\sqrt{-7}}2\times\f{a-b\sqrt{-7}}2=\f{a^2+7b^2}4$$
with $a,b\in\Z$ and $a\eq b\ (\mo\ 2)$.
Obviously $a$ and $b$ must be even (otherwise $a^2+7b^2\eq0\ (\mo\ 8)$), and $p=x^2+7y^2$ with $x=a/2$ and $y=b/2$.
\medskip

\proclaim{Conjecture 5.4} Let $p$ be an odd prime.
 Then
 $$\align&\sum_{k=0}^{p-1}\f{\bi{2k}k^2\bi{3k}k}{64^k}
 \\\eq&\cases
 x^2-2p\ (\mo\ p^2)&\t{if}\ (\f p{11})=1\ \&\ 4p=x^2+11y^2\ (x,y\in\Z),
 \\0\ (\mo\ p^2)&\t{if}\ (\f p{11})=-1.
 \endcases\endalign$$
 \endproclaim
 \Remark\ 5.3. It is well-known that the quadratic field $\Q(\sqrt{-11})$ has class number one
 and hence for any odd prime $p$ with $(\f p{11})=1$ we can write $4p=x^2+11y^2$ with $x,y\in\Z$.
 Concerning the parameters in the representation $4p=x^2+11y^2$, Jacobi
 obtained the following result (see, e.g., [BEW] and [HW]): If $p=11f+1$ is a prime and $4p=x^2+11y^2$ with $x,y\in\Z$ and $x\eq2\ (\mo\ 11)$, then
 $x\eq\bi{6f}{3f}\bi{3f}f/\bi{4f}{2f}\ (\mo\ p)$.

\proclaim{Conjecture 5.5} Let $p$ be an odd prime. If $p\eq1\pmod 4$ and $p=x^2+y^2$ with $x\eq1\pmod4$ and $y\eq0\pmod 2$,
then $$\sum_{k=0}^{p-1}\f{\bi{2k}k^2}{8^k}\eq\sum_{k=0}^{p-1}\f{\bi{2k}k^2}{(-16)^k}
\eq(-1)^{(p-1)/4}\l(2x-\f p{2x}\r)\ (\mo\ p^2)$$
and
$$\sum_{k=0}^{p-1}\f{\bi{2k}k^2}{32^k}\eq2x-\f p{2x}\ (\mo\ p^2).$$
If $p\eq 3\ (\mo\ 4)$ then
$$\sum_{k=0}^{p-1}\bi{p-1}k\f{\bi{2k}k^2}{(-8)^k}\eq\sum_{k=0}^{p-1}\f{\bi{2k}k^2}{32^k}\eq0\ (\mo\ p^2)$$
and
$$\sum_{k=0}^{p-1}\f{\bi{2k}k^2}{(-16)^k}\eq-\sum_{k=0}^{p-1}\f{\bi{2k}k^2}{8^k}\pmod{p^3}.$$
\endproclaim
\Remark\ 5.4. The author could prove all the congruences in Conjecture 5.5 modulo $p$.

\medskip

For more conjectures of the author on congruences related to central binomial coefficients, the reader may consult [Su4].

\Ack. The author is grateful to the referee for helpful comments.


\medskip

 \widestnumber\key{DPSW}

 \Refs

\ref\key BEW\by B. C. Berndt, R. J. Evans and K. S. Williams
\book Gauss and Jacobi Sums\publ John Wiley \& Sons, 1998\endref

\ref\key DPSW\by K. Driver, H. Prodinger, C. Schneider and J. Weideman
\paper Pad\'e approximations to the logarithm III: Alternative methods and additional results
\jour Ramanujan J.\vol 12\yr 2006\pages 299--314\endref

\ref\key F\by G. J. Fox\paper Congruences relating rational values of Bernoulli and Euler polynomials
\jour Fibonacci Quart.\vol 39\yr 2001\pages 50--57\endref

\ref\key G\by H. W. Gould\book Combinatorial Identities
\publ Morgantown Printing and Binding Co., 1972\endref


\ref\key GKP\by R. L. Graham, D. E. Knuth and O. Patashnik
 \book Concrete Mathematics\publ 2nd ed., Addison-Wesley, New York\yr 1994\endref

\ref\key GZ\by G. Grossman and A. Zeleke\paper On linear recurrence relations
\jour J. Concr. Appl. Math.\vol 1\yr 2003\pages 229--245\endref

\ref\key G-Z\by V. J. W. Guo and J. Zeng\paper Some congruences involving central $q$-binomial coefficients
\jour Adv. in Appl. Math.\vol 45\yr 2010\pages 303--316\endref

\ref\key HW\by R. H. Hudson and K. S. Williams\paper Binomial coefficients and Jacobi sums
\jour Trans. Amer. Math. Soc.\vol 281\yr 1984\pages 431--505\endref

\ref\key I\by T. Ishikawa\paper On Beukers' congruence\jour Kobe J. Math.\vol 6\yr 1989\pages 49--52\endref

\ref\key L\by E. Lehmer\paper On congruences involving Bernoulli numbers and the quotients
of Fermat and Wilson\jour Ann. of Math.\vol 39\yr 1938\pages 350--360\endref

\ref\key Ma\by R. Matsumoto\paper A collection of formulae for $\pi$
\jour on-line version is available from the website
{\tt http://www.pluto.ai.kyutech.ac.jp/plt/matumoto/pi\_small}
\endref

\ref\key Mo\by F. Morley\paper Note on the congruence
$2^{4n}\equiv(-1)^n(2n)!/(n!)^2$, where $2n+1$ is a prime\jour Ann.
Math. \vol 9\yr 1895\pages 168--170\endref

\ref\key PS\by H. Pan and Z. W. Sun\paper A combinatorial identity
with application to Catalan numbers \jour Discrete Math.\vol
306\yr 2006\pages 1921--1940\endref

\ref\key PWZ\by M. Petkov\v sek, H. S. Wilf and D. Zeilberger\book $A=B$ \publ A K Peters, Wellesley, 1996\endref

\ref\key Sp\by R. Sprugnoli\paper Sums of reciprocals of the central binomial coefficients
\jour Integers\vol 6\yr 2006\pages \#A27, 18pp (electronic)\endref

\ref\key St1\by R. P. Stanley\book Enumerative Combinatorics \publ
Vol. 1, Cambridge Univ. Press, Cambridge, 1999\endref

\ref\key St2\by R. P. Stanley\book Enumerative Combinatorics \publ
Vol. 2, Cambridge Univ. Press, Cambridge, 1999\endref

\ref\key SSZ\by N. Strauss, J. Shallit and D. Zagier
\paper Some strange $3$-adic identities\jour Amer. Math. Monthly
\vol 99\yr 1992\pages 66--69\endref

\ref\key S\by Z. H. Sun\paper Congruences concerning Bernoulli numbers and Bernoulli polynomials
\jour Discrete Appl. Math.\vol 105\yr 2000\pages 193--223\endref

\ref\key Su1\by Z. W. Sun\paper Binomial coefficients, Catalan numbers and Lucas quotients
\jour Sci. China Math.\vol 53\yr 2010\pages 2473--2488\endref

\ref\key Su2\by Z. W. Sun\paper Conjecture on a new series for $\pi^3$\jour A
Message to Number Theory List (sent on March 31, 2010). Available from the website
\newline \pages {\tt
http://listserv.nodak.edu/cgi-bin/wa.exe?A2=ind1003$\&$L=nmbrthry$\&$T=0$\&$P=1956}\endref

\ref\key Su3\by Z. W. Sun\paper $p$-adic valuations of some sums of multinomial coefficients
\jour Acta Arith. \vol 148\yr 2011\pages 63--76\endref

\ref\key Su4\by Z. W. Sun\paper Super congruences and Euler numbers
\jour Sci. China Math\pages to appear. {\tt http://arxiv.org/abs/1001.4453}\endref

\ref\key Su5\by Z. W. Sun\paper On Delannoy numbers and Schr\"oder numbers
\jour J. Number Theory, to appear. {\tt http://arxiv.org/abs/1009.2486}\endref

\ref\key ST1\by Z. W. Sun and R. Tauraso\paper New congruences for central binomial coefficients
\jour Adv. in Appl. Math.\vol 45\yr 2010\pages 125--148\endref

\ref\key ST2\by Z. W. Sun and R. Tauraso\paper On some new congruences for binomial coefficients
\jour Int. J. Number Theory\vol 7\yr 2011\pages 645--662\endref

\ref\key ZG\by X. Zhu and G. Grossman\paper On zeros of polynomial sequence
\jour J. Comput. Anal. Appl. \vol11\yr 2009\pages 140--158\endref

\endRefs
\enddocument